\documentclass{amsart}

\usepackage{graphicx}
\usepackage{amsfonts}
\usepackage{epsf}
\usepackage{amssymb}
\usepackage{amsmath}
\usepackage{amscd}

\newtheorem{theorem}{Theorem}[section]
\newtheorem{proposition}[theorem]{Proposition}

\newtheorem{remark}[theorem]{Remark} 
\newtheorem{example}[theorem]{Example}
\newtheorem{obstruction}[theorem]{Obstruction}

\newcommand{\cokerr}{\mbox{Coker} }

\newcommand{\s}{\mathfrak{s}}
\newcommand{\tk}{\mathfrak{t}}

\newcommand{\calc}{\mathcal{C}}
\newcommand{\calg}{\mathcal{G}}
\newcommand{\zz}{\mathbb{Z}}

\begin{document}

\title[Order in the concordance group]{Order in the concordance group and Heegaard Floer homology}
\author{Stanislav Jabuka, Swatee Naik} 
\keywords{Concordance, Heegaard Floer homology.}
\begin{abstract}{We use the Heegaard-Floer homology correction terms defined by 
Ozsv\'{a}th--Szab\'{o} to formulate a new obstruction for a knot to be of finite 
order in the smooth concordance group. This obstruction bears a formal 
resemblance to that of Casson and Gordon but is sensitive to the difference
between the smooth versus topological category. As an application we obtain 
new lower bounds for the concordance order
of small crossing knots.  }
\end{abstract}
\maketitle
\section{Introduction}
A knot $K$ in $S^3$ is called {\sl slice} if $(S^3, K) = \partial (B^4, D^2)$ where
$D^2$ is a 2-disk smoothly and properly embedded in the 4-ball $B^4$.
Knots $K_1$ and $K_2$ are called {\sl concordant} if 
$K_1\ \#\ \overline{K}_2$ is slice where $\overline{K}$
represents the mirror image of
$K$ with reversed string orientation. The set of concordance classes
of knots forms an Abelian group under the connected sum operation called the 
smooth concordance group and is denoted by $\mathcal{C}_1$. 
The order of $K$ in this group is the least
positive $n$ for which the connected sum of $n$ copies of $K$ is slice.

In this paper we use the correction terms for 3-manifolds stemming from 
Heegaard-Floer homology
to obstruct torsion in $\calc_1$. Specifically, we focus our attention on knots with 
10 or fewer crossings. Among these there are, as of this writing, 26 knots with unknown concordance 
order. Table \ref{tableofknots} below, courtesy of {\sl KnotInfo}\footnote{{\sl KnotInfo} is an online atlas of knots maintained 
by Charles Livingston. It can be found at {\tt http://www.indiana.edu/$\sim$knotinfo}.}, lists these knots along 
with lower bounds on their orders.
\begin{table}[ht]
\caption{} \label{tableofknots}
\begin{tabular}{|c|c||c|c||c|c|} \hline 
Knot $K$ & Order of $K$ & Knot $K$ & Order of $K$ &Knot $K$ & Order of $K$ \cr \hline \hline 
$8_{13}$ & $\ge 4$  & $10_{26}$ & $\ge 4$ &  $10_{102}$ & $\ge 4$ \cr \hline 
$8_{17}$ & $\ge 4$  & $10_{28}$ & $\ge 4$ & $10_{109}$ & $\ge 4$ \cr \hline 
$9_{14}$ & $ \ge 4$ & $10_{34}$ & $\ge 4$ & $10_{115}$ & $\ge 4$ \cr \hline
$9_{19}$ & $ \ge 4$ & $10_{58}$ & $\ge 4$ & $10_{118}$ & $\ge 4$ \cr \hline
$9_{30}$ & $ \ge 4$ & $10_{60}$ & $\ge 4$ & $10_{119}$ & $\ge 4$ \cr \hline
$9_{33}$ & $ \ge 4$ & $10_{79}$ & $\ge 4$ & $10_{135}$ & $\ge 4$ \cr \hline
$9_{44}$ & $ \ge 4$ & $10_{81}$ & $\ge 4$ & $10_{158}$ & $\ge 2$ \cr \hline
$10_{10}$ & $ \ge 4$ & $10_{88}$ & $\ge 4$ & $10_{164}$ & $\ge 4$ \cr \hline
$10_{13}$ & $ \ge 4$ & $10_{91}$ & $\ge 4$ &  & \cr \hline
\end{tabular}
\end{table}
\vskip2mm
The existing lower bounds from table \ref{tableofknots} 
have been determined by 
A. Tamulis \cite{tamulis}.

The obstruction we use (elucidated in 
section \ref{obstrsection} in detail) for a knot to be of order $n$ in $\calc_1$ applies, in principle, 
to all $n\ge 2$. However, computational complexity prevents us from checking the obstruction for $n>4$. 
Nonetheless, for $n=4$ the algorithm gives the following improvement on the above table:
\begin{theorem} \label{main}
The concordance order of any knot $K$ from the set of 14 knots
$$\left\{ 
\begin{array}{cc} 8_{13}, \, 9_{14}, \, 9_{19}, \, 9_{33}, \,  9_{44}, \,  10_{13},\, 10_{26}, \,
10_{28},  \cr
10_{34}, \, 10_{58},  \, 10_{60}, \, 10_{102},  \, 10_{119},  \, 10_{135}
\end{array}
\right\}$$
is at least 6. 
\end{theorem}
\begin{remark}
While the correction term obstruction (see section \ref{obstrsection}) 
bears a formal resemblance to that of Casson and Gordon
\cite{CG1, CG2}, it is nonetheless bound to be substantially different. 
While the Casson-Gordon obstruction
does not differentiate between the subtle distinction of smooth versus 
topological sliceness, our methods
are indeed sensitive to it. For example, the pretzel knot 
$P(7,-3,5)$ is topologically slice (and so all of
its Casson-Gordon obstructions vanish) but our methods can be applied
to show that its order in $\calc_1$ is infinite. 
Further examples of this type can be found in \cite{owensmanolescu}.
\end{remark}

The structure of $\calc_1$ is still rather poorly understood and virtually nothing is known about 
torsion in $\calc_1$. We briefly summarize the current state of understanding of $\calc_1$ and 
point out connections to the knots from theorem \ref{main}. 

There is a surjective homomorphism
$\Theta:\calc_1\rightarrow \calg$ (Levine, \cite{L1, L2}) 
from $\calc_1$ onto the {\sl algebraic concordance group} $\calg$, which consists of
Witt classes of Seifert forms under orthogonal sums. It is known that $\calg$ is  
isomorphic to the infinite direct sum
$$\calg \cong \zz^\infty \oplus\zz_2^\infty \oplus\zz_4^\infty. $$
The analogous homomorphism from odd dimensional 
concordance groups $\calc_{2n+1}$ (concordance classes of embeddings of $S^{2n+1}$ into $S^{2n+3}$) 
is an isomorphism for $n>1$ and it is injective onto an index 2 subgroup of $\calg$ when $n=1$. 
In the case of $\calc_1$ the kernel is nontrivial
as first proved by Casson and Gordon \cite{CG1, CG2}. In fact, the kernel of $\Theta$, referred to as the subgroup of 
{\sl algebraically slice knots}, is known to contain a subgroup isomorphic to 
$\mathbb{Z}^\infty \oplus \mathbb{Z}_2^\infty$ by  work of Jiang \cite{jiang} and Livingston \cite{Liv2}. 
All the knots in table \ref{tableofknots} 
map to order two elements in $\calg$ and are therefore of either 
infinite order or finite and even order in $\calc_1.$

\begin{remark}
Given the isomorphism 
of the  higher dimensional concordance groups $\calc_{2n+1}$
with the 
group $\zz^\infty \oplus\zz_2^\infty \oplus\zz_4^\infty, $
it is a reasonable guess to expect $\calc_1$ to exhibit 4-torsion 
elements (besides existing 2-torsion, see below) and perhaps no other finite
torsion. In view of this, while the bounds from theorem \ref{main} 
are only incrementally greater than those
from table \ref{tableofknots}, the increase of the bounds 
past order 4 is an important one.
\end{remark}

A negative amphicheiral knot, that is a knot which is isotopic to its
mirror image with reversed orientation, is clearly of concordance order
2. Other than this nothing is currently 
known about torsion in $\calc_1.$ 
In higher dimensions there are order 2 concordance classes
not represented by negative amphicheiral knots \cite{CM}. 
In dimension three it is unknown whether or not
the corresponding order 2 algebraic concordance 
classes have any concordance order 2 representatives.

Levine's set of invariants of algebraic concordance includes the
Tristram-Levine signatures, which being additive 
integral invariants, vanish for 
any knot representing a finite order algebraic concordance class.
Accordingly all of these invariants vanish for the knots from table \ref{tableofknots}. 

There are other, more subtle obstructions (\lq\lq subtler\rq\rq in the sense that they 
differentiate between the smooth and topological slice genus, a topic which we don't discuss here)
to a knot representing a torsion class in $\calc_1$:
\begin{align} \nonumber
\tau (K) & = \mbox{The Ozsv\' ath-Szab\' o $\tau$ invariant from Heegaard Floer homology \cite{peter11}. } \cr
s(K) & = \mbox{The Rasmussen invariant defined using Khovanov homology \cite{rasmussen1}.} \cr
\delta(K) & = \mbox{The $\delta$-invariant of Manolescu and Owens also defined using} \cr
&  \mbox{\phantom{iiii} Heegaard Floer homology \cite{owensmanolescu}.}
\end{align}
If either of these is non-vanishing, the knot $K$ is of infinite order in $\calc_1$. For the knots 
from table \ref{tableofknots} all three of these invariants are either known or are readily calculated 
and are all vanishing. 

Yet further information about the concordance order of knots comes from the following theorem 
proved in \cite{LN2} using the Casson-Gordon obstructions to sliceness.
\begin{theorem} \label{swateestheorem} Let $K$ be a knot in $S^3$ with 2-fold branched cover
$Y_K$.  If $H_1(Y_K;\zz) \cong \zz_{p^n} \oplus G$ with $p$ a prime congruent to
3  mod 4, $n$ odd and $p$ not dividing the order of $G$, then $K$ is of
infinite order in  $\calc_1$.\end{theorem}

This theorem gives a rather strong obstruction to being a torsion element in $\calc_1$, 
however, 
as is easy to check, none of the knots from table \ref{tableofknots} satisfy 
the hypothesis of theorem \ref{swateestheorem}. 

Additional obstructions to sliceness were obtained in \cite{KL}
using the twisted Alexander polynomials which relate to 
determinants of Casson-Gordon invariants. 
Using these Tamulis showed
in \cite{tamulis} that 
none of the knots from table \ref{tableofknots} have order 2 in $\calc_1.$

In summary, the knots from table \ref{tableofknots} are rather resilient 
to most of the known concordance invariants. It is in this sense that the use of Heegaard Floer 
homology in the proof of theorem \ref{main} 
is a significant new method, one which we hope will bear more fruit in the near future.

The remainder of the article is organized as follows. Section \ref{heegaard} reviews 
relevant parts of Heegaard Floer homology and reminds the reader of basic properties of the 
3-manifold correction terms $d(Y, \s)$. Section \ref{obstrsection} states the obstruction to being order $n$ in 
$\calc_1$ coming from the said correction terms. Section \ref{calcobstrterms}
explains how we calculated the correction terms for the double branched covers of the knots 
from table \ref{tableofknots}. Finally, section \ref{applying} explains how the results of theorem 
\ref{main} follow from our main obstruction. 

No originality is claimed on the material presented in sections \ref{heegaard}--\ref{calcobstrterms}. 
Our main obstruction \ref{obstruction} has been first observed by Ozsv\' ath and Szab\' o \cite{peter4} 
and been successfully used by other authors \cite{owensstrle, owensmanolescu}. Our contribution is 
the use of this obstruction 
to address long-standing questions about torsion in $\calc_1$.  

{\bf Acknowledgement } We would like to thank Peter Ozsv\'ath, Brendan Owens and Sa\v so Strle for 
many helpful discussions and for so generously sharing their expertise.  Special thanks are due to 
Charles Livingston for his feedback and for creating and maintaining the very helpful {\sl KnotInfo}. 
\section{Heegaard Floer homology} \label{heegaard}
This section serves as a reminder of some basic definitions and properties of the Heegaard 
Floer homology groups and the resulting correction terms for 3-manifolds. 
\subsection{The Heegaard Floer homology groups}
In their seminal papers \cite{peter1,peter2} Peter Ozsv\'ath and Zolt\'an Szab\'o introduced the 
Heegaard Floer homology groups $\widehat{HF}(Y,\s)$, 
$HF^\pm(Y,\s)$ and $HF^\infty(Y,\s)$ associated to a spin$^c$ 3-manifold 
$(Y, \s)$. These Abelian groups come equipped with a relative $\mathbb{Z}_d$--grading $gr$ where 
$$ d = \gcd \{ \langle c_1(\s),h\rangle \, | \, h\in H_2(Y;\mathbb{Z}) \}$$
In the case when $\s$ is torsion (by which we mean that $c_1(\s)$ is torsion) the relative 
$\mathbb{Z}$--grading $gr$ lifts to an absolute $\mathbb{Q}$--grading $\widetilde{gr}$. 
 
The various Heegaard Floer groups are related by means of long exact sequences. 
For example $HF^\pm (Y,\s) $ and $HF^\infty(Y,\s) $ fit into the sequence
\begin{equation} \label{LES}
 ... \rightarrow HF^-(Y,\s) \rightarrow HF^\infty (Y,\s) \stackrel{\pi}{\rightarrow} HF^+(Y,\s) 
\rightarrow HF^-(Y,\s) \rightarrow...
\end{equation}
If $\s$ is torsion then the maps in the above sequence preserve the absolute grading $\widetilde{gr}$ 
except the map 
$HF^+(Y,\s)\rightarrow HF^-(Y,\s)$ which drops degree by 1. 
\subsection{Cobordism induced maps}
The Heegaard Floer homology groups fit into a TQFT framework in the following sense: given a spin$^c$ 
4-manifold $(W,\tk)$ with $\partial W= -Y_1 \sqcup Y_2$ (where $-Y$ is $Y$ with its orientation reversed)
there are induced group homomorphisms 
$$ F^\circ _{W,\tk}  :HF^\circ (Y_1, \tk|_{Y_1}) \rightarrow HF^\circ (Y_2, \tk|_{Y_2}) $$
where $\circ$ stands for any of $\widehat{\phantom{n}}$, $+$, $-$, $\infty$. 
When $\tk|_{Y_1}$ and 
$\tk|_{Y_2}$ are both torsion the degree shift of the map $F_{W,\tk}^\circ$ is  
\begin{equation} \label{degshift}
\deg F_{W,\tk}^\circ := \widetilde{gr}(F_{W,\tk}^\circ (x)) - \widetilde{gr}(x) = 
\frac{(c_1(\tk))^2 -2e_W -3\sigma _W}{4}
\end{equation}
where $e_W$ and $\sigma _W$ are the Euler number and signature of $W$ respectively and 
$x\in HF^\circ (Y_1, \tk|_{Y_1})$ is any homogeneous element. Said differently, $F_{W,\tk}^\circ$ is a 
homogeneous map of degree $((c_1(\tk))^2 -2e_W -3\sigma _W)/4$.

\begin{proposition} [Ozsv\'ath-Szab\'o, \cite{peter4}] \label{negdefcobor}
When $b_2^+(W) = 0$ the homomorphism $F_{W,\tk}^\infty$ is an isomorphism for all spin$^c$-structures 
$\tk$ on $W$. 
\end{proposition}

The exact sequence \eqref{LES} is functorial under cobordism induced maps in the sense that one 
obtains the commutative diagram (with exact rows):
\begin{equation} \label{commdiag}
\begin{CD} 
 @>>> HF^-(Y_1,\s_1) @>>> HF^\infty(Y_1,\s_1) @>\pi>> HF^+(Y_1,\s_1) @>>>  \\
@.  @VF_{W,\tk}^-VV  @VF_{W,\tk}^\infty VV @VF_{W,\tk}^+VV  \\
@>>> HF^-(Y_2,\s_2) @>>> HF^\infty(Y_2,\s_2) @>\pi>> HF^+(Y_2,\s_2) @>>>  \\
\end{CD}
\end{equation}
In the above diagram $\s_i$ stands for $\tk|_{Y_i}$.
\subsection{The correction terms for 3-manifolds}
Let $Y$ be a rational homology sphere and let $\s \in Spin^c(Y)$ be a spin$^c$-structure on $Y$. 
The {\sl correction term} $d(Y,\s)$ is defined to be 
$$ d(Y,\s) = \min \{ \widetilde{gr}(\pi(x))  \, | \, x \in HF^\infty (Y,\s) \}$$
where $\pi :HF^\infty (Y,\s) \rightarrow HF^+(Y,\s)$ is the map from the exact sequence \eqref{LES}.

\begin{example} \label{thesphere}
Consider $S^3$ with its unique spin-structure $\s_0$. Recall from \cite{peter1} that 
$HF^\infty(S^3,\s_0) \cong \mathbb{Z}[U,U^{-1}]$ and 
$HF^+(S^3,\s_0) \cong \mathbb{Z}[U^{-1}]$. 
The absolute grading on both groups is 
specified by $\widetilde{gr}(U^k) = -2k$ and the map $\pi : HF^\infty(S^3,\s_0) \rightarrow HF^+(S^3,\s_0)$
is the obvious quotient map 
$$ \mathbb{Z}[U,U^{-1}] \rightarrow \frac{\mathbb{Z}[U,U^{-1}]}{U \,\mathbb{Z}[U]} \cong \mathbb{Z}[U^{-1}]$$
Thus $\pi$ is surjective and therefore $d(S^3,\s_0)$ is the lowest grading 
in $HF^+(S^3,\s_0)$ which in turn is given by 
\begin{equation} \label{s3}
d(S^3,\s_0) = \widetilde{gr}(U^0) = 0
\end{equation}
\vskip1mm
\end{example}

The correction terms enjoy a number of nice properties. Given $\s \in Spin^c(Y)$ let $\overline{\s}$ be
the conjugate spin$^c$-structure.  Then  
\begin{align} \label{corecprop}
d(Y,\overline{\s}) & = d(Y,\s) \cr
d(-Y,\s) & = -d(Y,\s)\cr
d(Y_1\#Y_2,\s_1\#\s_2) & = d(Y_1,\s_1) + d(Y_2,\s_2)
\end{align}
%
\subsection{Correction terms for 3-manifolds bounding rational homology 4-balls}
Consider now two rational homology 3-spheres $Y_1$ and $Y_2$ equipped with spin$^c$-structures  
$\s_i\in Spin^c(Y_i)$. Consider furthermore a 
negative definite cobordism $(W,\tk)$ from $(Y_1,\s_1)$ to $(Y_2,\s_2)$ (i.e a 4-manifold $W$ 
with $\partial W = -Y_1\sqcup Y_2$, $\tk|_{Y_i} = \s_i$ and $b_2^+(W) = 0$).  
Let $x_2\in HF^\infty (Y_2,\s_2)$ 
be an element with $\widetilde{gr}(\pi(x_2)) = d(Y_2,\s_2)$ where $\pi$ is the map 
from \eqref{LES}. According to proposition 
\ref{negdefcobor} the homomorphism $F_{W,\tk}^\infty:HF^\infty (Y_1,\s_1) \rightarrow HF^\infty (Y_2,\s_2)$ 
is an isomorphism. Let $x_1 \in HF^\infty (Y_1,\s_1)$ be the unique preimage of $x_2$ under this map. 
The degree-shift formula \eqref{degshift}  and the commutative diagram \eqref{commdiag} show that 
$$ \widetilde{gr}(\pi(x_2)) - \widetilde{gr}(\pi(x_1)) = \frac{(c_1(\tk))^2 - 2e_W -3\sigma _W}{4}$$
Since $d(Y_1,\s_1) \le \widetilde{gr}(\pi(x_1))$ by definition and $d(Y_2,\s_2) = \widetilde{gr}(\pi(x_2))$
by choice of $x_2$, the above equality becomes the inequality 
\begin{equation} \label{theinequality}
d(Y_1,\s_1) \le d(Y_2,\s_2) - \frac{(c_1(\tk))^2 - 2e_W -3\sigma _W}{4}
\end{equation}

Let us now turn to the special case when $Y_2=S^3$ and $W$ has the rational homology of a punctured 4-ball. 
Then $e_W = 0=\sigma _W$ and $(c_1(\tk))^2=0$ for all $\tk \in Spin^c(W)$. The above inequality along with 
example \ref{thesphere} then shows  
that $d(Y_1,\s_1)\le 0 $ for all spin$^c$-structures $\s_1 \in Spin^c(Y_1)$ which lie in the image of the 
restriction map $Spin^c(W) \rightarrow Spin^c(Y_1)$. Reversing the orientation on $W$ and applying 
\eqref{theinequality} once more shows that $d(Y_1,\s_1)\ge 0 $ also. Therefore 
$$ d(Y_1,\s_1) = 0  \quad \quad \forall \, \s_1 \in \mbox{Im} [Spin^c(W) \rightarrow Spin^c(Y_1)]$$    
If we fill in the $S^3$ boundary component of $W$ with a 4-ball, we see that $Y=Y_1$ bounds a rational homology 
ball $X = W\cup_{S^3} B^4$. It is well known (and follows easily from the universal coefficient theorem and 
the exact sequence of the pair $(X,Y)$) that such a 3-manifold has second cohomology of square order, say 
$|H^2(Y;\mathbb{Z})|=n^2$, and that the order of the image $H^2(X;\mathbb{Z})\rightarrow H^2(Y;\mathbb{Z})$ 
is $n$. After suitable affine identifications of $Spin^c(X)\cong H^2(X;\mathbb{Z})$ and $Spin^c(Y) 
\cong H^2(Y;\mathbb{Z})$ the restriction map $Spin^c(X)\rightarrow Spin^c(Y)$ corresponds precisely to the 
restriction induced map $H^2(X;\mathbb{Z}) \rightarrow H^2(Y;\mathbb{Z})$. We summarize our discussion in the 
following 
\begin{theorem} \label{thevanishing}
Let $Y$ be a rational homology 3-sphere which bounds a rational homology 4-ball $X$. 
Then $|H^2(Y;\mathbb{Z})|=n^2$ for some $n$ and there is a subgroup $\mathcal{P}$ of 
$H^2(Y;\mathbb{Z})$ of order $n$ such that 
$$d(Y,\s)=0  \quad \quad \forall \, \s \in \mathcal{P}$$
under a suitable identification $Spin^c(Y) \cong H^2(Y;\mathbb{Z})$.   
\end{theorem}
%
%
\section{The sliceness obstruction} \label{obstrsection}
Let $K$ be a knot in $S^3$ and let $Y_K$ be the double branched cover of $S^3$ branched along $K$. The order 
of the second cohomology of $Y_K$ is given by 
$$ |H^2(Y_K;\mathbb{Z})| =  |\det (K)| = |\Delta _K(-1)|$$
where $\Delta_K(t)$ is the 
Alexander polynomial of $K$. 

If $K$ is slice with slice disk $D^2 \hookrightarrow B^4$ we let $X_K$ be the double branched cover of $B^4$ 
branched along $D^2$. The manifold $X_K$ is a rational homology ball with boundary $\partial X_K = Y_K$. 
Thus according to theorem \ref{thevanishing} we must have $|\det(K)| = n^2$ for some integer $n$ and 
$d(Y_K, \s) = 0$ for all $\s$ in some subgroup $\mathcal{P}$ of $H^2(Y_K;\mathbb{Z})$ of order $n$. 
As sample calculations show, this turns out to be a rather strong obstruction to the sliceness of $K$. 

To apply this algorithm to the question of the order of a knot $K$ in $\mathcal{C}_1$  consider the knot 
$K' = \#^{2m} K$, the $2m$-fold connected sum of $K$ with itself. 

If $K$ is of order $2m$ then $K'$ is slice and 
the above algorithm asserts the vanishing of $d(Y_{K'},\s')$ for spin$^c$-structures $\s'$ from 
some (affine) subgroup 
$\mathcal{O}$ of $Spin^c(Y_{K'})$ of order $|\det(K)|^m$. Recall that 
$$ Y_{K_1 \# K_2} \cong Y_{K_1} \# Y_{K_2}\quad \quad Spin^c(Y_1 \# Y_2 ) \cong Spin^c(Y_1) \times Spin^c(Y_2) $$
Thus a spin$^c$-structure $\s' \in Spin^c(Y_{K'})$ corresponds to a collection of 
$2m$ spin$^c$-structures $\s' = (\s_1,...,\s_{2m})$ with $\s_i \in Spin^c(Y_K)$.  Furthermore \eqref{corecprop} 
implies that for such an $\s'$ the correction term $d(Y_{K'},\s')$ is given by 
$$ d(Y_{K'},\s') = d(Y_K,\s_1) + ... + d(Y_K,\s_{2m})$$

To summarize we obtain the following 
\begin{obstruction} \label{obstruction}
If $K$ is of order $2m$ in the smooth knot concordance group $\mathcal{C}_1$  
there exists a subgroup $\mathcal{O}$ of $\left( H^2(Y;\mathbb{Z})\right)^{2m} \cong (Spin^c(Y_K))^{\times 2m}$ 
of order $|\det(K)|^m$ with 
\begin{equation} \label{vanishingsum}
d(Y_K,\s_1) + ... + d(Y_K,\s_{2m}) = 0 \quad \quad \forall (\s_1,...,\s_{2m}) \in \mathcal{O}
\end{equation}
In the above $Y_K$ is the double branched cover of $S^3$ branched along $K$.
\end{obstruction}

One drawback of this obstruction algorithm is that there is a priori no way of knowing what the 
group $\mathcal{O}$ might
be in the case of a concrete knot $K$. We are thus forced to consider all subgroups $\mathcal{O}$ of 
$H^2(Y_K;\mathbb{Z})^{\times 2m}$ of order $|\det (K)|^m$ and hope that for none of them relation \eqref{vanishingsum} holds. 
If this is the case, $K$ cannot be of order $2m$. 

To use this obstruction for the knots $K$ from table \ref{tableofknots} one needs to calculate the 
correction terms $d(Y_K,\s)$ for all $\s \in Spin^c(Y_K)$. 
In the next section we do this by distinguishing a number of cases. 
\section{Calculating obstruction terms} \label{calcobstrterms}
\subsection{2-bridge knots} \label{2bridge}
Some of the knots from table \ref{tableofknots} are 2-bridge knots and so their double branched covers
are lens spaces. 
The correction terms for lens spaces have been calculated 
by Ozsv\'ath and Szab\'o in \cite{peter4} and follow the recursive formula 
$$d(-L(p,q),i) = \left(\frac{pq-(2i+1-p-q)^2}{4pq} \right) - d(-L(q,r),j)$$
where $r$ and $j$ are the mod $q$ reductions of $p$ and $i$ respectively. Here $i$ is an integer 
$0\le i < p+q$ whose mod $p$ reduction represents the spin$^c$-structure $[i]\in \mathbb{Z}_p \cong
Spin^c(-L(p,q))$. 

The knots from table \ref{tableofknots} whose double branched covers $Y_K$ are lens spaces are
\vskip1mm
\begin{center}
\begin{table}[ht]
\caption{} \label{table2}
\begin{tabular}{|c|c||c|c|} \hline 
Knot $K$ & $Y_K$ & Knot $K$ & $Y_K$  \cr \hline \hline 
$8_{13}$ & $L(29,11)$ & $10_{13}$ & $L(53,22)$ \cr \hline
$9_{14}$ & $L(37,14)$ & $10_{26}$ & $L(61,17)$ \cr \hline
$9_{19}$ & $L(41,16)$ & $10_{28}$ & $L(53,19)$ \cr \hline
$10_{10}$ & $L(45,17)$ & $10_{34}$ & $L(37,13)$ \cr \hline
\end{tabular}
\end{table}
\end{center}
\vskip1mm
For example, the correction terms of $Y_{8_{13}}$ thus obtained are 
\begin{align} \nonumber
\left\{ -\frac{2}{29},-\frac{18}{29},\frac{8}{29},\frac{18}{29},\frac{12}{29}, -\frac{10}{29},\frac{10}{29},
\frac{14}{29},\frac{2}{29}, -\frac{26}{29},-\frac{12}{29},-\frac{14}{29},-\frac{32}{29},-\frac{8}{29}, 0,
\right. \cr
\quad \left. -\frac{8}{29}, 
-\frac{32}{29},-\frac{14}{29},-\frac{12}{29},-\frac{26}{29}, 
\frac{2}{29},\frac{14}{29},\frac{10}{29}, -\frac{10}{29},\frac{12}{29},\frac{18}{29},
\frac{8}{29},-\frac{18}{29},-\frac{2}{29},\right\} 
\end{align}
\subsection{Alternating knots} \label{alternating}
When $K$ is a knot which possesses an alternating projection $D$  the correction terms of $Y_K$ 
can be calculated from the Goeritz matrix $G$ associated to $D$. The details of this have been worked out 
by Ozsv\'ath and Szab\'o in \cite{peter17} and we summarize them here for the benefit of the reader. 

\begin{figure}[htb!] 
\centering
\includegraphics[width=3cm]{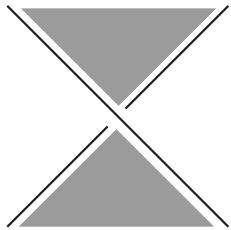}
\caption{The coloring conventions near a crossing.}  \label{pic1}
\end{figure}

Let $D$ be an alternating projection of a knot $K$. We color the regions of $D$ black and white 
according to the convention from figure \ref{pic1}, to obtain
a checkerboard pattern. 

\begin{figure}[htb!] 
\centering
\includegraphics[width=10cm]{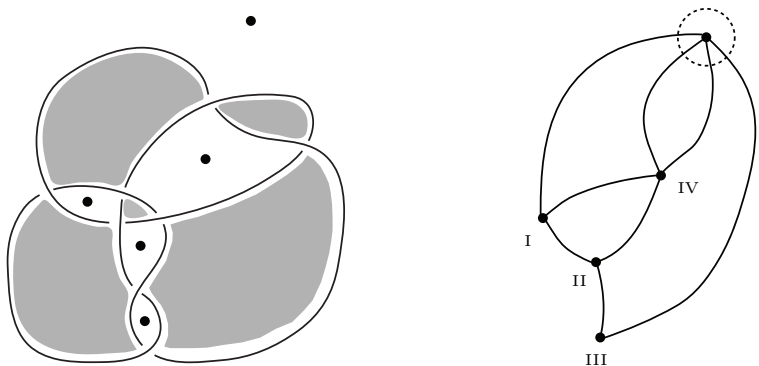}
\put(-65,0){\tiny III}
\put(-70,30){\tiny II}
\put(-88,45){\tiny I}
\put(-30,65){\tiny IV}
\caption{The checkerboard diagram and graph associated to the knot $8_{17}$. In the construction of the 
Goeritz matrix we drop the vertex enclosed by the dotted circle.}  \label{pic3}
\end{figure}

From such a pattern we extract a graph in the following way: The vertices of the graph 
are in bijection with the white regions in the diagram (including the unbounded region if it 
happens to be white). There is an edge between two vertices for each touching 
point of their corresponding white regions. 
Figure \ref{pic3}  shows the checkerboard diagram and the associated graph for the knot $8_{17}$. 

From the graph we now 
extract a matrix - the Goeritz matrix of the projection $D$. Pick and discard 
one of the vertices of the graph (the vertex enclosed in a dotted circle in figure \ref{pic3}). Give the remaining
vertices an arbitrary ordering. The Goeritz matrix $G=[g_{ij}]$ has the entries
$$ g_{ij} = \left\{ 
\begin{array}{lc}
\mbox{Number of edges between the $i$-th and $j$-th vertex} & \quad ; \quad i\ne j \cr
-1\cdot \mbox{Valence of the $i$-th vertex} & \quad ; \quad i= j
\end{array}
\right.
$$
For example, the Goeritz matrix associated to the projection of $8_{17}$ from figure \ref{pic3} with the ordering 
of the vertices as indicated is 
$$ G_{8_{17}} = \left[
\begin{array}{rrrr}
-3 & 1 & 0 & 1 \cr
1 & -3 & 1 & 1 \cr
0 & 1 & -2 & 0 \cr
1 & 1 & 0 & -4
\end{array}
\right]
$$
Finally, from the Goeritz matrix it is now a matter of arithmetic to extract the correction terms 
$d(Y_K,\s)$: Consider $G:V\otimes V \rightarrow \mathbb{Z}$ as a 
negative definite bilinear quadratic form where $V=\mathbb{Z}^{\ell}$ if $G$ is of dimension $\ell \times \ell$. 
Let $g:V\rightarrow V^*$ and $G^*:V^*\otimes V^*\rightarrow \mathbb{Q}$ 
be the obvious maps induced by $G$, namely 
$$ g(v)= G(v,\cdot)  \quad \quad \mbox{ and } \quad \quad G^*(G(v,\cdot),G(w,\cdot)) = G(v,w) $$
Let $M_g:\cokerr(g) \rightarrow \mathbb{Q}$ be
$$ M_g(\xi) = \frac{1}{4} \left( \max_{\{ v\in V^* \, |\, [v]=\xi\} } G^*(v_0 + 2v, v_0+2v) + \mbox{rk}(V) \right)$$
where $v_0$ is any characteristic vector of $G$ (i.e. any vector $v_0\in V^*$ 
with $v_0(w) \equiv G(w,w)\,\, (\mbox{mod}\, \,  2)$ for all $w\in V$). It is shown in \cite{peter17} 
that there is an 
isomorphism $ \varphi : \cokerr (g) \rightarrow H^2(Y_K;\mathbb{Z})$ such that 
$$ d(Y_K,\varphi(\xi)) = M_g(\xi)$$
for some affine identification of $H^2(Y_K;\mathbb{Z})$ with $Spin^c(Y_K)$.

For example, the correction terms for $Y_{8_{17}}$ calculated this way are 
\begin{align} \nonumber
\left\{ -\frac{20}{37},-\frac{32}{37},\frac{18}{37},-\frac{18}{37},\frac{8}{37},\frac{22}{37},\frac{24}{37},
\frac{14}{37},-\frac{8}{37},\frac{32}{37},-\frac{14}{37},\frac{2}{37},\frac{6}{37},\right. \cr
\left. -\frac{2}{37},-\frac{22}{37},\frac{20}{37},-\frac{24}{37},-\frac{6}{37},0,-\frac{6}{37},-\frac{24}{37},
\frac{20}{37},-\frac{22}{37},-\frac{2}{37},\frac{6}{37}, \right. \cr
\left. \frac{2}{37},
-\frac{14}{37},\frac{32}{37},-\frac{8}{37},\frac{14}{37},\frac{24}{37},\frac{22}{37},\frac{8}{37},
-\frac{18}{37},\frac{18}{37},-\frac{32}{37},-\frac{20}{37} \right\}
\end{align}

Of the knots from table \ref{tableofknots} which do not appear in 
table \ref{table2}, the 
ones which have alternating projections are listed in table \ref{table3}. 

\begin{center}
\begin{table} \caption{} \label{table3}
\begin{tabular}{|c|} \hline 
Knots having alternating projections \cr\hline \hline
$8_{17}$, $9_{30}$, $9_{33}$, $10_{58}$, $10_{60}$, $10_{79}$, $10_{81}$ \cr \hline 
$10_{88}$, $10_{91}$, $10_{102}$, $10_{109}$, $10_{115}$, $10_{118}$, $10_{119}$ \cr \hline
\end{tabular}
\end{table}
\end{center}
\vskip2mm
%

\subsection{The remaining cases}
Sections \ref{2bridge} and \ref{alternating} allow for a calculation of the correction terms of $Y_K$ 
for most knots $K$
from table \ref{tableofknots}. The knots from that table which do not fall into either category are 
\begin{equation} \label{specknots}
9_{44} \quad 10_{135} \quad 10_{158} \quad 10_{164} 
\end{equation}
and these require special attention. All four of these knots however \lq\lq resemble\rq\rq 
alternating knots sufficiently 
so that a calculation of their correction terms can be done by using the Goeritz matrix again. 

\begin{figure}[htb!] 
\centering
\includegraphics[width=6cm]{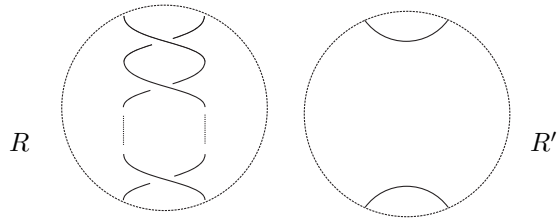}
\put(-190,25){$R$}
\put(7,25){$R'$}
\caption{The left-handed 
orientation of the twists in the region $R$ is shown on the left. The region $R'$ on the right
is used to replace $R$ in forming $L$ from $K$.}  \label{pic4}
\end{figure}

The following algorithm has been described in \cite{peter17}, see also \cite{owensmanolescu}. 
Suppose that $K$ is a knot with a knot projection $D$ which outside some region $R$ is alternating and inside 
$R$ consists of $k$ left-handed half-twists of two parallel strands, see figure \ref{pic4}.  
Such projections can be found for all four knots from \eqref{specknots}. 

\begin{figure}[htb!] 
\centering
\includegraphics[width=8cm]{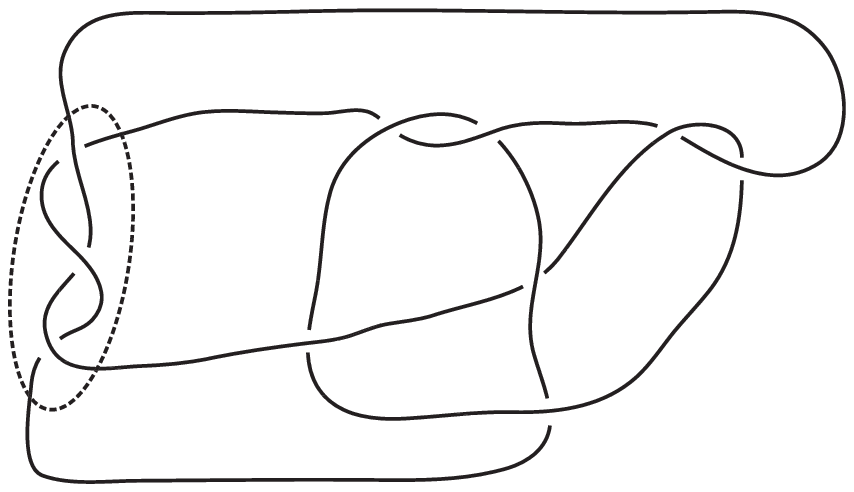}
\put(-240,60){$R$}
\put(-30,50){$10_{158}$}
\caption{The region $R$ is the portion of this projection of $10_{158}$ inside the dotted oval.}  \label{pic5}
\end{figure}

\begin{figure}[htb!] 
\centering
\includegraphics[width=10cm]{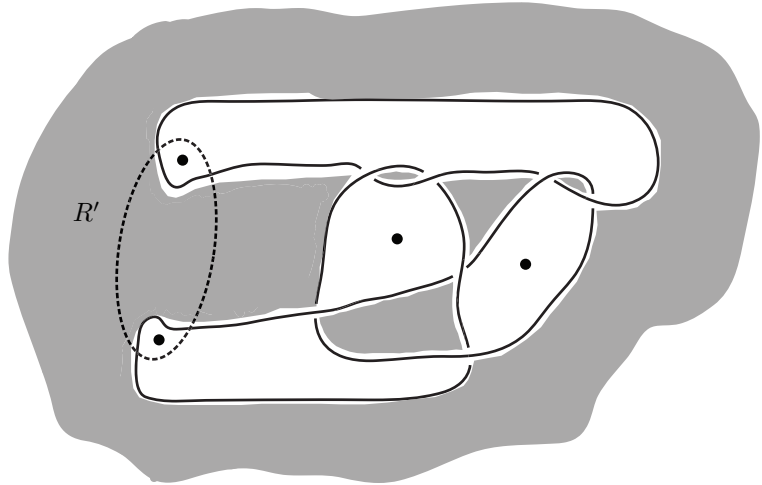}
\put(-260,100){$R'$}
\caption{The link $L$ in case of $K=10_{158}$.}  \label{pic6}
\end{figure}

Let $L$ be the 2-component link obtained from $K$ by replacing $R$ with $R'$.
For example, figure \ref{pic5} shows a knot projection of 
$10_{158}$ with the marked region $R$, figure \ref{pic6} depicts the corresponding link $L$.

%
%


%
%

There is a 
restriction we impose: the vertex 
from the checkerboard pattern for $L$ that we drop in the 
computation of the Goeritz matrix of $L$,
should always be one of the vertices from 
the region $R$. In figure \ref{pic6}, two such vertices are
indicated. 

Let $\tilde G$ be the Goeritz matrix of $L$ and let $G$ be the 
matrix obtained from $\tilde G$ as 
\begin{equation} \label{newgoeritz}
G = \left[ 
\begin{array}{rrrr|r}
 & & & & 0 \cr
 & \tilde G& & & \vdots \cr
 & & & & 0 \cr
 & & & & 1 \cr \hline
 0& \hdots & 0& 1& -k \cr
\end{array} 
\right]
\end{equation}
where $k$ is the number of negative half-twists in the region $R$.

When $K$ is any of $9_{44}$, $10_{135}$, $10_{158}$ or $10_{164}$
the correction terms of $Y_{K}$ are calculated from $G$ 
in the way described in section \ref{alternating} for alternating knots.  

For example, the Goeritz matrix $\tilde G$ of the link $L$ (figure \ref{pic6}) and its associated 
matrix $G$ for the knot $K=10_{158}$ (figure \ref{pic5}) are 
$$ \tilde G = \left[ 
\begin{array}{rrr}
-4 & 1 & 2  \cr
1 & -4 & 2  \cr
2 & 2 & -4  
\end{array}
\right]
\quad \quad \quad \quad 
G = \left[ 
\begin{array}{rrrr}
-4 & 1 & 2 & 0 \cr
1 & -4 & 2 & 0 \cr
2 & 2 & -4 & 1 \cr 
0 & 0 & 1 & -3  
\end{array}
\right]
$$
leading to the correction terms 
\begin{align} \nonumber
&\left\{ -\frac{2}{45},-\frac{2}{5},\frac{8}{9},-\frac{8}{45},\frac{2}{5},\frac{28}{45},\frac{22}{45},0,
-\frac{38}{45},-\frac{2}{45},\frac{2}{5},\frac{22}{45},
\frac{2}{9},-\frac{2}{5},\frac{28}{45},\right. \cr
& -\frac{32}{45},-\frac{2}{5},-\frac{4}{9},-\frac{38}{45},\frac{2}{5},-\frac{32}{45},-\frac{8}{45},0,-\frac{8}{45},
-\frac{32}{45},\frac{2}{5},-\frac{38}{45},-\frac{4}{9},-\frac{2}{5},-\frac{32}{45},\cr
& \quad \quad \quad 
\quad \quad \left. \frac{28}{45},-\frac{2}{5},\frac{2}{9},\frac{22}{45},
\frac{2}{5},-\frac{2}{45},-\frac{38}{45},0,\frac{22}{45},\frac{28}{45},\frac{2}{5},-\frac{8}{45},
\frac{8}{9},-\frac{2}{5},-\frac{2}{45} \right\} 
\end{align}
%
%
%
\section{Applying the obstruction} \label{applying}
Given a knot $K$, the obstruction \ref{obstruction} implies that if $K$ has order $2m$ in 
$\calc_1$ then there is a 
subgroup $\mathcal{O}$ of $H^2(Y_K;\mathbb{Z})^{\times 2m}$ of order $|\det(K)|^m$ for which all corresponding
correction terms vanish. To check the obstruction for a concrete knot $K$ one needs to:
\begin{enumerate}
\item Calculate all correction terms of $Y_K$. 
\item Find all subgroups $\mathcal{O}$ of $H^2(Y_K;\mathbb{Z})^{\times 2m}$ of order $|\det(K)|^m$.
\item Check that 
$$ d(Y_K,\s_1)+...+d(Y_K,\s_m) = 0 \quad \quad \quad \forall \, (\s_1,...,\s_m) \in \mathcal{O}$$
\end{enumerate}
We have written a {\sc Mathematica} script which performs each of the 3 steps above. Computationally the
most demanding part by far is step 2. In fact, our computational resources only allowed us to use $m=2$ (and thus 
test for 4-torsion in $\calc_1$) and even in that case we were forced to use a weaker version of obstruction 
\ref{obstruction}:
\begin{obstruction} \label{secondary}
If $|\det (K)| = p$ or $|\det (K)| = p\cdot q$ where $p\ne q$ are primes, then if $K$ is of order 4 there 
exists a subgroup $\widetilde{\mathcal{O}}$ of 
$H^2(Y_{K'};\mathbb{Z})$ isomorphic to  $\mathbb{Z}_{|\det (K)|}$ 
with 
$$ d(Y_{K},\s_1)+d(Y_{K},\s_2)+d(Y_{K},\s_3)+d(Y_{K},\s_4)= 0 $$
for all $(\s_1,\s_2,\s_3,\s_4) \in \widetilde{\mathcal{O}}$. Here $K'$ denotes $\#^4K$. 
\end{obstruction}   
This is a direct consequence of obstruction \ref{obstruction}. 
The results of theorem \ref{main} follow from 
our {\sc Mathematica} implementation of obstruction 5.1. 

An easy check reveals that all knots from 
table \ref{tableofknots} satisfy the hypothesis of obstruction \ref{secondary} except $10_{10}$, $10_{158}$ and
$10_{164}$. Each of these 3 knots has determinant $45$. In these cases 
it still follows from obstruction \ref{obstruction}  
that there is a subgroup $\widetilde{\mathcal{O}}$ of $H^2(Y_{K'};\mathbb{Z})$ of order 45
whose associated correction terms vanish. However, unlike in obstruction \ref{secondary}, 
there are now 2 possibilities for the isomorphism type of 
$\widetilde{\mathcal{O}}$, namely $\mathbb{Z}_{45}$ and $\mathbb{Z}_3 \oplus \mathbb{Z}_{15}$. 
While for
each of the knots $10_{10}$, $10_{158}$ and $10_{164}$ no group of the former 
type (with vanishing correction terms) exists, there are groups of the latter
type and so no conclusions can be drawn.




\end{document}